\documentclass{amsart}
\usepackage{amssymb}
\usepackage{mathrsfs}
\usepackage{graphicx}
\usepackage{comment}
\usepackage{rotating}

\newcommand{\be}{\begin{equation}}
\newcommand{\ee}{\end{equation}}
\newcommand{\ba}{\begin{align}}
\newcommand{\ea}{\end{align}}

\newtheorem{theorem}{Theorem}[section]

\newtheorem{lemma}{Lemma}[section]

{\begin{list}{}{%
\settowidth{\labelwidth}{\textsf{{\it #1.}}}%
\setlength{\labelsep}{4mm}%
\setlength{\leftmargin}{\labelwidth}%
\addtolength{\leftmargin}{\labelsep}%
}}%
{\end{list}}

\def\om{\omega}
\def\beq{\begin{equation}}\def\enq{\end{equation}}

{\begin{list}{}{%
\settowidth{\labelwidth}{\textsf{{\it #1.}}}%
\setlength{\labelsep}{2mm}%
\setlength{\leftmargin}{\labelwidth}%
\addtolength{\leftmargin}{\labelsep}%
\addtolength{\leftmargin}{4mm}%
\setlength{\itemsep}{6pt}%
\setlength{\listparindent}{0pt}%
\setlength{\topsep}{3pt}%
}}%

\title[Integer group determinants ]{The integer group determinants for $GA(1,q)$}

\author[A. Ostergaard]{Andrew Ostergaard}
\address{Department of Mathematical Sciences\\
University of Southern Indiana\\
Evansville, IN 47712, USA}
\email{a.ostergaard@usi.edu}

\author[C. Pinner]{Chris Pinner}
\address{ Department of Mathematics\\
         Kansas State University\\
         Manhattan, KS 66506, USA}
\email{pinner@math.ksu.edu}


\date{\today}
\begin{document}

\begin{abstract}
We show that the integer group determinants for the general affine group of degree one, $GA(1,q)$ with $q=p^k$ a prime power, take the form $D=AB^{q-1},$ where $A$ is a $\mathbb Z_{q-1}$ integer group determinant and $B\equiv A \bmod q$. This generalizes the result for $k=1$. 

When $2^k-1$ is a Mersenne prime we show that this condition  is both necessary and sufficient for $GA(1,2^k).$ The same is true for $GA(1,9)$ and $GA(1,27)$, but interestingly not  for $GA(1,16)$.


\end{abstract}

\maketitle

\section{Introduction}
The integer group determinant problem for a finite group $G=\{g_1,\ldots,g_n\}$  is to determine which integers are achieved as a group determinant all of whose entries are integer. That is,  thinking of this as a function on elements of the group ring $\mathbb Z [G]$ and defining
$$ D\left(  \sum_{g\in G} a_g g\right)=\det \left(a_{g_i g_j^{-1}} \right), $$
$i, j=1,\ldots,n$ indexing rows and columns respectively, we want to describe
$$ S(G)=\left\{ D\left(  \sum_{g\in G} a_g g\right)\; :\; a_{g_1},\ldots ,a_{g_n}\in \mathbb Z\right\}. $$
For cyclic groups $G=\mathbb Z_m$, Newman \cite{Newman1} 
and  Laquer  \cite{Laquer} showed that
\be \label{cyclicachieve} \{ n\in \mathbb Z \; : \; \gcd(m,n)=1 \text{ or } m^2\mid n\} \subseteq S(\mathbb Z_m), \ee
together with  the divisibilty restrictions on primes $r$ dividing a $\mathbb Z_m$ integer group determinant $D$ and $m$:
\be \label{cyclicdiv} r^{\alpha} \parallel m, \;\;r\mid D \; \Rightarrow \; r^{\alpha+1}\mid D. \ee
When $G=\mathbb Z_p$ or $\mathbb Z_{2p},$  $p$ a prime, the restrictions  \eqref{cyclicdiv} are if and only if. Even for small groups, the general problem can quickly become difficult,  although a complete description of $S(G)$ has been obtained for all groups with $|G|<20$ (see \cite{Bishnu18} for references).

Here we are interested in the  general affine groups of degree one
$$ GA(1,q)=\left\{ \begin{pmatrix} a & b \\ 0 & 1 \end{pmatrix}\; : \; a\in \mathbb F_q^*,\; b\in \mathbb F_q\right\} $$
where $q=p^k$, $p$ prime. This should be thought of as the semi-direct product
$\mathbb F_q \rtimes_{\theta} \mathbb F_q^*$ where the automorphism $\theta(x)$ corresponds to multiplication by $x$ in $\mathbb F_q$.

When $k=1$ it was shown in \cite{GA1P} that the integer group determinants take the form
\be \label{k=1}  D=AB^{p-1}, \; \text{ $A$ is a $\mathbb Z_{p-1}$ integer determinant, } B\equiv A \bmod p. \ee

A similar structure  holds for $GA(1,q)$ for prime powers $q=p^k$.

\begin{theorem} \label{main} Let $p$ be a prime and $q=p^k,$ then the integer group determinants for $GA(1,q)$ take the form
\be \label{genk}  D=AB^{q-1},  \; \text{ $A$ is a $\mathbb Z_{q-1}$ integer determinant, } B\equiv A \bmod q. \ee
\end{theorem}
It is tempting to ask whether the restriction \eqref{genk} is sufficient as well as necessary.
From \cite{GA1P} we know that this is true for $GA(1,p)$ when $p=3,5,7,11,23$.
In Theorem \ref{Mersenne}  below, we show that it is also holds when $q=2^k$ with $q-1$ a Mersenne prime, and
in Sections \ref{q=9} and \ref{q=27} we verify it for $GA(1,3^2)$ and $GA(1,3^3)$.
However (see Section \ref{q=16}) it turns out not to be true for $GA(1,16)$.


\section{A group presentation for $GA(1,q)$}

We shall suppose that $r$ is a generator of $\mathbb F_q^*$ with minimal polynomial 
$$ f(x)=x^k-a_{k-1}x^{k-1} - \cdots -a_0,\;\; \;\;a_0,\ldots ,a_{k-1}\in \{0,1,\ldots, p-1\}, $$
over $\mathbb F_p$. We write
$$ X:=\begin{pmatrix} r & 0 \\ 0 & 1 \end{pmatrix}, \;\; \;\;Y_j=\begin{pmatrix} 1 & r^{j-1} \\ 0 & 1 \end{pmatrix}, \; j=1,\ldots , k. $$
Observe that any $a\in \mathbb F_q^*$ and $b\in \mathbb F_q$ can be written
$$ a=r^{\ell},\;\; b=b_1+b_2r+\cdots + b_{k} r^{k-1}, $$
for some $0\leq \ell < q-1$, $0\leq b_1,\ldots ,b_{k} <p,$ and hence 
$$ \begin{pmatrix} a & b \\ 0 & 1 \end{pmatrix}= Y_1^{b_1}\cdots Y_{k}^{b_{k}} X^{\ell}. $$
Hence we shall think of our group ring elements as taking the form
\be \label{polyform} F(X;Y_1,\ldots ,Y_{k})= \sum_{j=0}^{q-2} f_j(Y_1,\ldots , Y_{k})X^j \ee
where the 
$$f_j(Y_1,\ldots , Y_{k})=\sum_{i_1=0}^{p-1}\cdots \sum_{i_{k}=0}^{p-1} a_j(i_1,\ldots ,i_{k}) Y_1^{i_1}\cdots Y_{k}^{i_{k}}$$ 
are polynomials in $\mathbb Z[Y_1,\ldots ,Y_{k}]$.
Since $XY_jX^{-1}=\begin{pmatrix} 1 & r^{j} \\ 0 & 1 \end{pmatrix} $ we have 
\be \label{gprels1} X^{q-1}=Y_1^p=\cdots = Y_{k}^p=1, \quad Y_iY_j=Y_jY_i \text{ for any $1\leq i,j\leq k$,}\ee
and
\be \label{gprels2} XY_j=Y_{j+1}X, \; \text{ for $j=1,\ldots ,k-1$},\;\;\;\; XY_{k}=Y_1^{a_0}\cdots Y_{k}^{a_{k-1}} X. \ee
That is, we can think of $GA(1,q)$ as being generated by $X,Y_1,\ldots, Y_{k},$ 
subject to the group relations \eqref{gprels1} and \eqref{gprels2}. Of course everything in the group can be written in terms of the two generators $X$ and $Y_1$.

\section{Formula for the group determinant}
We recall Frobenius' formula \cite{Frob} giving  a factorization of the group determinant using a complete set $\hat{G}$ of irreducible representations for $G$:
$$ D\left(  \sum_{g\in G} a_g g\right)=\prod_{\rho\in \hat{G}} \det \left(\sum_{g\in G}  a_g \rho(g)\right)^{\deg \rho}. $$
In this case we will have $q$ representations (see for example \cite{Conrad1}); 
from 
$$\begin{pmatrix} x & y \\ 0 & 1 \end{pmatrix}\begin{pmatrix} a & b \\ 0 & 1 \end{pmatrix}\begin{pmatrix} x & y \\ 0 & 1 \end{pmatrix}^{-1}=\begin{pmatrix} a & bx-y(a-1) \\ 0 & 1 \end{pmatrix}$$
 just observe that there are the $q$ conjugacy classes:
$$ \left\{ \begin{pmatrix} a & b \\ 0 & 1 \end{pmatrix} \; : \; b\in \mathbb F_q\right\} \text{ for $a\in \mathbb F_q^* \setminus \{1\}$}, \quad \left\{ \begin{pmatrix} 1 & b \\ 0 & 1 \end{pmatrix}\; :\; b\in \mathbb F_q^*\right\}, \quad \left\{\begin{pmatrix} 1 & 0 \\ 0 & 1 \end{pmatrix}\right\}.$$
As shown below, there are $(q-1)$ multiplicative characters (the degree one representations), leaving one degree $(q-1)$ representation, $q(q-1)=(q-1)\cdot 1 + 1\cdot (q-1)^2$.

Suppose that $\chi$ is a multiplicative character on $G$.
From \eqref{gprels1} we see that
$$ \chi(X)^{q-1}=1, \;\;\; \chi (Y_1)^p=1,$$
and from \eqref{gprels2}
$$\chi(Y_1)=\chi(Y_2)=\cdots =\chi(Y_{k}),\;\;\;\chi(Y_1)^{1-a_0-\cdots -a_{k-1}}=1. $$
Now plainly $f(x)$ does not have a factor $x-1$, $f(1)\not\equiv 0 \bmod p$ and $\chi (Y_1)^p=\chi(Y_1)^{f(1)}=1$ forces $\chi(Y_1)=1$.

That is we have $q-1$ characters; those with $\chi(X)$ being any $(q-1)$st root of unity and
$\chi(Y_1)=\cdots =\chi(Y_{k})=1.$ These  contribute the integer
$$ A = \prod_{x^{q-1}=1} F(x;1,\cdots, 1) $$
to the group determinant for \eqref{polyform}. Note, this is the $\mathbb Z_{q-1}$
determinant for $F(x;1,\ldots ,1)$.

Now we attempt to find a  degree $(q-1)$ representation $\rho,$ where $\rho(X)$ is the basic circulant shift matrix, $P,$ and for $j=1,\ldots ,k$, the $\rho(Y_j)$ are diagonal matrices with first entry $y_j$ (note these entries will all be $p$th roots of unity);
\be \label{defP}  \rho(X) =P:=\begin{pmatrix} 0  & 1 & 0 &  \cdots & 0 \\ 0 & 0 &1 & \cdots &  0 \\
  & \vdots  & & \ddots & \vdots \\
0 & 0 & 0 & \cdots &  1\\
1 & 0 & 0 & \cdots &  0 \end{pmatrix}, \hspace{2ex} \rho(Y_j) = \begin{pmatrix} y_j & & & \\
 & \ddots &  & \\
 & & \ddots & \\
 & &  & \ddots\end{pmatrix}. \ee

Observe that if we write
$$XY_1^{u_1}\cdots Y_{k}^{u_{k}}X^{-1}=Y_1^{a_0u_{k}}Y_2^{u_1+a_{1}u_{k}} \cdots Y_{k}^{u_{k-1}+a_{k-1}u_{k}}$$
as $Y_1^{v_1}\cdots Y_{k}^{v_{k}},$ then
$$ \begin{pmatrix} v_1\\ \vdots \\ v_{k}\end{pmatrix} = M\begin{pmatrix} u_1\\ \vdots \\ u_{k}\end{pmatrix},\quad M:=\begin{pmatrix} 0 & 0 & 0  & \cdots & 0 &  a_0\\
1 & 0 & 0  &   & 0 & a_1 \\
0 & 1 &  0&  & 0 &  a_2 \\
\vdots & & & & \vdots  & \vdots \\
0 & 0 & 0  & \cdots &1 & a_{k-1} \end{pmatrix}.$$
Note,  for a diagonal matrix $T$ the effect of $\rho(X)T\rho(X^{-1})$ is to circulantly
shift the diagonal sequence back one place.
Hence the entry in the $j$th row of $\rho(Y_1^{u_1}\cdots Y_{k}^{u_{k}})$ will be the first entry in $\rho\left( X^{j-1}Y_1^{u_1}\cdots  Y_{k}^{u_{k}}X^{-(j-1)}\right)$, that is $y_1^{v_1}\cdots y_{k}^{v_{k}}$  where
$$ \begin{pmatrix} v_1\\ \vdots \\ v_
{k}\end{pmatrix} = M^{j-1} \begin{pmatrix} u_1\\ \vdots \\ u_{k}\end{pmatrix}, \quad j=1,\ldots ,q-1.$$
In particular the $j$th entry in $\rho(Y_1)$ will be $y_1^{u_1}\cdots y_{k}^{u_{k}}$
where
\be \label{geny} \begin{pmatrix} u_1 \\  \vdots \\ u_{k} \end{pmatrix}= M^{j-1}\begin{pmatrix} 1 \\ 0 \\ \vdots \\ 0\end{pmatrix} \bmod p, \quad j=1,\ldots , q-1. \ee
It is easily seen, expanding $xI-M$ along its last column,  that the characteristic equation for $M$ is $f(x)$ which has distinct roots. Hence, writing
 $$M=K\begin{pmatrix} \lambda_1 & & & \\
 & \lambda_2 &  & \\
 & & \ddots & \\
 & & & \lambda_{q-1} \end{pmatrix}K^{-1} $$ 
with distinct non-zero $\lambda_i$ in $\mathbb F_q^*,$ then plainly the $\lambda_i^{q-1}=1$ and $M^{q-1}=I$ and the sequence of exponents generated by \eqref{geny} has period $q-1$ and never hits $(0,\ldots ,0)$. 
 Moreover the roots of $f(x)$ all have order $q-1$ (they must be roots of the same cyclotomic factor of $x^{q-1}-1$), the eigenvalues of $M^{\ell}$ are
not 1 for any $1\leq \ell <q-1,$ and the  sequence of exponents cannot have order less than $q-1$. In particular the diagonal entries of $\rho(Y_1)$ will run through all the $y_1^{u_1}\cdots y_{k}^{u_{k}},$ with $ 0\leq u_1,\ldots, u_{k}<p$ except for $u_1=\ldots =u_{k}=0$.

Likewise we could obtain the exponents  for the entries in  $\rho(Y_2),\ldots, \rho(Y_{k})$, from 
$$M^{j-1}\begin{pmatrix} 0\\1\\ \vdots \\ 0\end{pmatrix},\ldots , M^{j-1}\begin{pmatrix} 0\\0\\ \vdots \\ 1\end{pmatrix};$$
 that is, the $\ell$th column of $M^{j-1}$ gives the exponents on $y_1,\ldots ,y_{k}$ for the element in the  $j$th row in $\rho(Y_{\ell })$. But notice the diagonals for $\rho(Y_2),\rho(Y_3),\ldots, $  will actually just be a circulant shift one, two, \ldots, places of the diagonal entries for $\rho(Y_1)$  already generated by \eqref{geny}.

Thus, for any choice of $p$th roots of unity $y_1,\ldots ,y_{k},$ a $\rho$ constructed
this way will satisfy the group relations and give us a representation.
We write
$\vec{y}_1=(y_1,\ldots ,y_{k}),$ and
$\vec{y}_2=(y_2,\ldots,y_{k},y_1^{a_0}\cdots y_{k}^{a_{k-1}}),\ldots ,\vec{y}_{q-1},$
for the corresponding monomials in the $j$th row, $j=1,\ldots ,q-1$, of $\rho(Y_1),\ldots ,\rho(Y_{k})$.
Notice that 
$\rho( Y_1^{\alpha_1}\cdots Y_{k}^{\alpha_{k}})$ will have $\vec{y}_j(1)^{\alpha_1}\cdots \vec{y}_j(k)^{\alpha_{k}},$  down the diagonal
for $j=1,\ldots ,q-1$. Similarly 
$\rho (f(Y_1,\ldots ,Y_{k}))$ will be a diagonal matrix with  $f(\vec{y}_j)$ in the $j$th row and $\rho (f(Y_1,\ldots ,Y_{k})X^t)$ will be this matrix with a circulant shift $t$ places to the right. Hence $\rho$ will contribute the term
\be \label{defB} B=\det \begin{pmatrix} f_0(\vec{y}_1) & f_1(\vec{y}_1) & \cdots & f_{q-2}(\vec{y}_1) \\ f_{q-2}(\vec{y_2}) & f_0(\vec{y}_2) & \cdots & f_{q-3}(\vec{y}_2) \\ \vdots &  \vdots &  & \vdots \\  f_1(\vec{y}_{q-1}) & f_2(\vec{y}_{q-1}) & \cdots & f_0(\vec{y}_{q-1}) \end{pmatrix}. \ee
See the $q=9$ example in Section \ref{q=9}.

Notice any choice of $p$th roots $y_1,\ldots ,y_{k}$ gives a representation, but if
we take  all of them 1, then  the $\vec{y}_j(t)$ will all be 1 and the entries in \eqref{defB} all take the 
form $f_j(1,\ldots ,1).$ This determinant is the circulant  determinant corresponding to the $\mathbb Z_{q-1}$ determinant of $F(X;1,\ldots,1)$ and the determinant \eqref{defB} just gives us $A$ again (i.e. the contribution coming from the one dimensional representations).

If $y_{\ell}\neq 1$ say, then taking $F(X;Y_1,\ldots ,Y_{k})=\Phi_{p}(Y_{\ell}),$ $\Phi_p(x)=1+\cdots + x^{p-1}$ denoting the $p$th cyclotomic polynomial, the first row of \eqref{defB} will be all zeros
and the determinant $B=0$, while  $F(X;1,\ldots 1)=p$ and $\chi(F(X;Y_1,\ldots ,Y_{k}))=p$ for any of the characters $\chi$.  Hence  $B\neq p^{q-1}$ cannot come from a product of $(q-1)$ characters, but from the degree $q-1$ representation.

Notice, we could also write the matrix for $B$ in terms of additive characters. If
\be \label{addchar} \psi(r^{i-1})=y_i, \;\; 1\leq i \leq k, \quad \psi(b_1+b_2r+\cdots +b_kr^{k-1}) =y_1^{b_1}\cdots y_k^{b_k}, \ee
then $\psi(x)$ will be an additive character of $\mathbb F_q$ for any choice of $p$th roots of $y_1,\ldots ,y_k$ (non-trivial if the $y_i$ are not all 1). The cycle of variables 
down the matrix for $B$ is then $\psi(r^i)$, $i=0,...,q-2$. Since $r^i$, $i=0,\ldots,q-2,$
generates all of $\mathbb F_q^*,$ we will cycle once through all $q-1$ monomials $y_1^{\alpha_1}\cdots y_k^{\alpha_k}$, $0\leq \alpha_i<p$, $\alpha_i$ not all 0.  The entry in the $(i+1)$st
row and $(j+1)$st column of the matrix for $B$ becomes
$$ f_{\ell}\left(\psi(r^i),\psi(r^{i+1}),\ldots ,\psi(r^{i+k})\right),\;\; \ell\equiv j-i\text{ mod}{ \:(q-1)}, \quad i,j=0,\ldots ,q-2. $$
In this form it is easy to see that $B$ is independent of the choice of $y_1,\ldots ,y_k$
(not all 1). Recall (see, for example \cite[Theorem 5.7]{LidlNied}) that the non-trivial additive characters all take the form $\psi (cx)$ for some $c$ in $\mathbb F_p^*$.
Hence using a different additive character, $\psi(r^u x)$ for some $1\leq u\leq q-2$, merely causes a circulant shift of $u$ places in the rows and the columns, and no change in the determinant $B$. In particular the algebraic integer $B$ is unchanged by the automorphisms $\om \mapsto \om^j$, $1\leq j<p$ of $\mathbb Q (\om),$ $\om:=\exp(2\pi i/p)$, and so is an integer.

\section{Another way to write this}\label{alternative}
We order the non constant monomials $y_1^{\alpha_1}\cdots y_k^{\alpha_k}$, $0\leq \alpha_i<p$ not all 0, using 
$$ \mathbf{y}_j=y_1^{u_1(j)}\cdots y_k^{u_k(j)},\quad \quad \begin{pmatrix} u_1(j) \\ \vdots \\ u_k(j)\end{pmatrix}  = M^j \begin{pmatrix} 1 \\ 0\\ \vdots \\ 0 \end{pmatrix} \pmod p, \quad j=0,\ldots ,q-2.$$
In terms of the additive character \eqref{addchar} this is just
$$ \mathbf{y}_j = \psi(r^j),\;\; j=0,\ldots ,q-2. $$
If we write 
$ F(x;y_1,\ldots ,y_k) = \sum_{j=0}^{q-2} f_j(y_1,\ldots ,y_k) x^j $
in the form
$$ f_j(y_1,\ldots ,y_k) =c_j + \sum_{n=0}^{q-2} c_{jn}(1 -\mathbf{ y}_n),$$
 then $F(x;1,\ldots ,1)=\sum_{j=0}^{q-2} c_jx^j$ and $A$ is its circulant determinant
\be \label{defP}  A=\det {\mathcal A},\quad A= \sum_{j=0}^{q-2} c_jP^j,\ee
where $P$ is the usual circulant shift matrix \eqref{defP}. Writing $S$ for the diagonal matrix with entries
$$1-\mathbf{y}_0, \ldots ,1-\mathbf{y}_{q-2}, $$
(when $k=1$ these are just the $1-y^{g^j\pmod p}$, $j=0,\ldots ,q-2$, where $g$ is  a primitive root mod $p$), and using $P^nSP^{-n}$ to circulantly shift the diagonal down $n$ places and $P^j$ to circulantly shift $j$ places to the right to the $x^j$ position, we have $D=AB^{q-1}$ with 
$$ B=\det{\mathcal B},\quad\quad {\mathcal B}= \mathcal{A} + \sum_{j=0}^{q-2}\sum_{n=0}^{q-2} c_{jn}P^nSP^{j-n}$$
after reducing mod $\langle y_1^p-1,\ldots ,y_k^p-1,\Phi_p(y_1)\cdots \Phi_p(y_k)\rangle$ or  putting in any choice of $p$th roots of unity for $y_1,\ldots ,y_k$ (other than $y_1=\cdots =y_k=1$ which, of course, gives $A$).

\section{Proof of Theorem \ref{main}}
The determinant \eqref{defB} is a polynomial in $y_1,\ldots ,y_{k}$. Reducing mod $y_j^p-1$ for $j=1,\ldots ,k$ we can assume it takes the form
$$ g(y_1,\ldots ,y_{k}) = \sum_{i_1=0}^{p-1}\cdots \sum_{i_{k}=0}^{p-1}a(i_1,\ldots ,i_{k}) y_1^{i_1}\cdots y_{k}^{i_{k}}\in \mathbb Z[y_1,\ldots ,y_{k}]. $$
Observing that
$$ \sum_{y^p=1} y^{t}= \begin{cases} 0 & \text{ if $t=1,\ldots ,p-1$,}\\ p & \text{ if $k=0$, } \end{cases} $$
we plainly have (the term $y_1=\cdots =y_{k}=1$ giving $A$ and the others all $B$)
\begin{align*}
A + (p^k-1) B  & = \sum_{y_1^p=1} \cdots \sum_{y_{k}^p=1} g(y_1,\ldots ,y_{k})\\ & = p^k a(0,\ldots ,0), 
\end{align*}
and $B\equiv A \bmod p^k$ as claimed.

\section{Achieving some values}
Recall, one can achieve all integers coprime to 
$m$ and all multiples of $m^2$ as $\mathbb Z_m$ integer determinants.
The following  two Theorems can be thought 
of as the  $GA(1,q)$ counterpart of  these classical cyclic results \eqref{cyclicachieve}. 

\begin{theorem} \label{coprime} We achieve all values of the form
$$ D=A B^{q-1},  \quad \gcd(A,q-1)=1,\quad B\equiv A \bmod q, $$
as  $GA(1,q)$ integer group determinants.

\end{theorem}

\begin{proof} We write $A=\ell +\lambda (q-1)$, with $1\leq \ell < q-1$, $\gcd(\ell,q-1)=1$
and observe that
$$ g(x)=1+x+\cdots + x^{\ell-1} + \lambda (1+x+ \cdots + x^{q-2}) $$
has $\mathbb Z_{q-1}$ determinant $\ell + \lambda (q-1)=A$.

We take
$$ F(x;y_1,\cdots ,y_k)=1+x+\cdots +x^{\ell-1} + t(y_1,\ldots ,y_k)(1+x+\cdots + x^{q-2}). $$
Observe that the determinant to obtain $B$ for $F$ resembles the determinant giving $B$ 
for $1+\cdots +x^{\ell-1}$, call this $B_G$,  except that we add $t(y_1,\ldots ,y_k)$ to each entry in the first
row (and their shifts under $M$  in the subsequent rows); subtracting the last column from the first $q-2$ columns we see that $B$ will be linear in $t(y_1,\ldots, y_k)$ 
with coefficient $a$ obtained by replacing the 
the first row in the determinant for $B_G$ by all 1's. Likewise for the subsequent rows with the same coefficient $a$.
Letting  $\vec{y}_1=(y_1,\ldots ,y_k)$, $\vec{y}_2$,\ldots ,$\vec{y}_{q-1},$
denote the cycle of variables in the rows of the matrix defining $B$, we get
$$ B=B_G + a( t(\vec{y}_1)+t(\vec{y}_2)+\cdots + t(\vec{y}_{q-1}) ). $$
Now  $B_G$ is just $\ell$, the $\mathbb Z_{q-1}$ determinant for $1+\cdots +x^{\ell-1}$. Taking $t(\vec{y}_1)=\lambda$ we get $A$, the $\mathbb Z_{q-1}$ determinant of $g(x)$ as above, and so $a=1$.

That is
$$ B=\ell +  t(\vec{y}_1)+t(\vec{y}_2)+\cdots + t(\vec{y}_{q-1}). $$
where the values $\vec{y}_1,\ldots ,\vec{y}_{q-1}$ run through all $(q-1)$ possible $k$-tuples of $p$th roots of unity
except for $(1,\ldots ,1)$. Hence taking 
$$t(y_1,\ldots ,y_k)=\lambda + m(1-y_1) $$
we get 
$$ B=\ell +\lambda (q-1) + mq=A+mq. \qedhere$$

\end{proof}

\begin{theorem} \label{square} We achieve all values of the form
$$ D=A B^{q-1},  \quad (q-1)^2 \mid A,\quad B\equiv A \bmod q, $$
as  $GA(1,q)$ integer group determinants.

\end{theorem}

\begin{proof} We take
$$ F(X; Y_1,\ldots ,Y_{k})= (1-Y_1X)+ t(Y_1,\ldots ,Y_{k})(1+X+\cdots + X^{q-2}),$$
so that $A,$ the $\mathbb Z_{q-1}$ determinant of $F(X;1,\ldots ,1)=(1-X)+t(1,\ldots ,1)(X^{q-1}-1)/(X-1),$ is
$$ A=t(1,\ldots ,1)(q-1)^2. $$
Writing $C_1=y_1,C_2,\ldots ,C_{q-1}$ for the first elements of $\vec{y}_1,\ldots ,\vec{y}_{q-1},$
$$ B=B_0+ \sum_{j=1}^{q-1}t(\vec{y}_j)\alpha(\vec{y}_j), $$
where $B_0$ is the value of $B$ for $1-Y_1X$
$$ B_0=\det \begin{pmatrix} 1 & -C_1 & 0 & \cdots & 0 \\
0 & 1 & -C_2 & \cdots & 0 \\
\vdots &  & & & \vdots \\ -C_{q-1} & 0 & 0 & \cdots & 1 \end{pmatrix} =1-C_1\cdots C_{q-1}=0, $$
and
\begin{align*} \alpha(y_1,\ldots ,y_{k}) & =\det \begin{pmatrix} 1 & 1 & 1 & \cdots & 1 \\
0 & 1 & -C_2 & \cdots & 0 \\
\vdots &  & & & \vdots \\ -C_{q-1} & 0 & 0 & \cdots & 1 \end{pmatrix}\\
 & =1+C_{q-1}+C_{q-1}C_{q-2}+ \cdots + C_{q-1}C_{q-2}\cdots C_2.
\end{align*}
We observe that we cannot have a product of  $t<q-1$ consecutive $C_i$ equalling one mod $\langle y_1^p-1,\ldots ,y_{k}^p-1\rangle. $ To have $C_j=y_1^{u_1}\cdots y_{k}^{u_{k}}$ with $C_j\cdots C_{j+t}=1$ mod $\langle y_1^p-1,\ldots ,y_{k}^p-1\rangle $, we must have 
$$ \left( I + M_T +  \cdots + M_T^{t-1}\right) \begin{pmatrix} u_1\\ \vdots \\ u_{k}\end{pmatrix} = \begin{pmatrix} 0 \\ \vdots \\ 0\end{pmatrix} \pmod p.$$
But we know that the eigenvalues $\lambda$ of $M_T$ are of order $(q-1)$ in $\mathbb F_q$, and so $1+\lambda +\cdots +\lambda^{t-1}\neq 0$ for $t<q-1$. 

Hence the terms in $\alpha(y_1,\ldots ,y_{k})$ will give $(q-1)$ of the possible $y_1^{u_1}\cdots y_{k}^{u_{k}}$, $0\leq u_1,\ldots ,u_{k}<p$ and if $y_1^{U_1}\cdots y_{k}^{U_{k}}$ is the missing monomial then
\begin{align*} \alpha(y_1,\ldots ,y_{k}) & =-y_1^{U_1}\cdots y_{k}^{U_{k}}, \\ \left(1-y_1^{p-U_1}\cdots y_{k}^{p-U_{k}}\right)\alpha(y_1,\ldots ,y_{k}) &=1-y_1^{U_1}\cdots y_{k}^{U_{k}}, \end{align*}
and
$$ t(x)=c+ \ell  \left(1-y_1^{p-U_1}\cdots y_{k}^{p-U_{k}}\right) $$
gives $A=c(q-1)^2$, $ B=c+\ell q$.
\end{proof}

\section{ The $GA(1,2^k)$ when $2^k-1$ is prime. }

When $2^k-1=r$ is a Mersenne prime, the integer group determinants for $\mathbb Z_{q-1}$ are exactly the $A$ with $\gcd(A,r)=1$ or $r^2\mid A$, and for these $A$ we can
achieve any the $B\equiv A$ mod $q$ from
Theorems \ref{coprime} and \ref{square}.

\begin{theorem} \label{Mersenne} If $2^k-1=r$ is prime, then the integer group determinants for $GA(1,2^k)$ are exactly the
$$ D=A B^{q-1}, \quad  r\nmid A  \text{ or } r^2\mid A,\quad B\equiv A \bmod q. $$
\end{theorem}

\section{The group $GA(1,9)$}\label{q=9}
We show that \eqref{genk} is also both a necessary and sufficient condition for $GA(1,9)$.

\begin{theorem} The integer group determinants for $GA(1,9)$ are exactly the
$$ D=A B^{8}, \quad  2\nmid A  \text{ or } 2^5 \mid A,\quad B\equiv A \bmod q. $$
\end{theorem}

\begin{proof} Writing $Y,Z$ instead of  $Y_0,Y_1$, an element in the group ring takes the form
$$ F(X;Y,Z)= \sum_{j=0}^{7} f_j(Y,Z)X^j $$
and
$$ A= \prod_{x^8=1} F(x;1,1). $$
We take $r=1+i$, $f(x)=x^2-2x-1,$ giving $M=\begin{pmatrix} 0 & 1 \\ 1 & 2\end{pmatrix}$. For $j=0,\ldots ,7,$  the
$ M^j \begin{pmatrix} 1 \\ 0 \end{pmatrix}$
run through
$  \begin{pmatrix} 1 \\ 0 \end{pmatrix},  \begin{pmatrix} 0 \\ 1  \end{pmatrix},
 \begin{pmatrix} 1 \\ 2 \end{pmatrix},  \begin{pmatrix} 2 \\ 2 \end{pmatrix},  \begin{pmatrix} 2 \\ 0 \end{pmatrix}, \begin{pmatrix} 0 \\ 2 \end{pmatrix}, \begin{pmatrix} 2 \\ 1 \end{pmatrix}, \begin{pmatrix} 1 \\ 1 \end{pmatrix} \pmod 3.$
That is,  down the rows of $B$ the $y$ and $z$ cycle through
$$ y, z,yz^2,y^2z^2,y^2,z^2,y^2z,yz,\text{ and }  z,yz^2,y^2z^2,y^2,z^2,y^2z,yz,y,$$
and $B$ is the determinant of
\begin{tiny}
$$  \left[ \begin{matrix} f_0(y,z) & f_1(y,z) & f_2(y,z) & f_3(y,z) & f_4(y,z) & f_5(y,z) & f_6(y,z) & f_7(y,z) \\
f_7(z,yz^2) & f_0(z,yz^2) &  f_1(z,yz^2) & f_2(z,yz^2) & f_3(z,yz^2) & f_4(z,yz^2) & f_5(z,yz^2)   & f_6(z,yz^2) \\
f_6(yz^2,y^2z^2) & f_7(yz^2,y^2z^2) & f_0(yz^2,y^2z^2) & f_1(yz^2,y^2z^2) & f_2(yz^2,y^2z^2) & f_3(yz^2,y^2z^2) & f_4(yz^2,y^2z^2)   & f_5(yz^2,y^2z^2) \\
f_5(y^2z^2,y^2) & f_6(y^2z^2,y^2) & f_7(y^2z^2,y^2) & f_0(y^2z^2,y^2) & f_1(y^2z^2,y^2) & f_2(y^2z^2,y^2) & f_3(y^2z^2,y^2) & f_4(y^2z^2,y^2)\\
f_4(y^2,z^2) & f_5(y^2,z^2) & f_6(y^2,z^2) & f_7(y^2,z^2) & f_0(y^2,z^2) & f_1(y^2,z^2) & f_2(y^2,z^2) & f_3(y^2,z^2)\\
f_3(z^2,y^2z) & f_4(z^2,y^2z) & f_5(z^2,y^2z) & f_6(z^2,y^2z) & f_7(z^2,y^2z) & f_0(z^2,y^2z) & f_1(z^2,y^2z) & f_2(z^2,y^2z)\\
f_2(y^2z,yz) & f_3(y^2z,yz) & f_4(y^2z,yz) & f_5(y^2z,yz) & f_6(y^2z,yz) & f_7(y^2z,yz) & f_0(y^2z,yz) & f_1(y^2z,yz)\\
f_1(yz,y) & f_2(yz,y) & f_3(yz,y) & f_4(yz,y) & f_5(yz,y) & f_6(yz,y) & f_7(yz,y) & f_0(yz,y)
\end{matrix}\right]
$$
\end{tiny}
mod $\langle y^3-1, z^3-1\rangle$, for any choice of cube roots of unity $y,z \neq 1,1$.

From \cite{Norbert2} the $\mathbb Z_{8}$ determinants are exactly the $A$ with $A$ odd
or $2^5\mid A$. From Theorems \ref{coprime} and \ref{square} we can obtain any
$D=AB^8$ with $B\equiv A \pmod 9$ for the  $A$ odd or a $A$ a  multiple of $2^6$. It remains to get all
the $B\equiv A \pmod 9$ when $A$ is an odd multiples of $2^5$. 

From \cite{Norbert2} we can achieve 32 as a $\mathbb Z_8$ determinant using $1+x^2+x^3+x^4$. We take
$$ F(X;Y,Z) = 1+YX^2+ X^3+X^4 + (c+ b(1-Y^2Z^2))(1+X+\cdots + X^7). $$
Then $F(x;1,1)=1+x^2+x^3+x^4 + c(1-x^8)/(1-x)$ has $\mathbb Z_8$ determinant
$$ A = 2^3 (2^2+8c)=2^5(1+2c). $$
Writing $B_0$ for the $B$ value for $1+YX^2+X^3+X^4$, we have
\begin{align*} B_0 &  =\det \begin{pmatrix}1&0&y&1&1&0&0&0\\0 & 1 &  0 & z & 1 & 1 & 0 &0\\
0 & 0 & 1 & 0 & yz^2 & 1 &1 &0 \\
0 & 0 & 0 & 1 & 0 & y^2z^2 &1 & 1 \\
1 & 0 & 0 & 0 & 1 &0 & y^2 & 1 \\ 
1 & 1 & 0 &0 &0 & 1 & 0 & z^2\\
y^2z & 1 & 1 & 0 & 0 & 0 & 1 &0 \\ 
0 & yz & 1 & 1 & 0 & 0 & 0 & 1 \end{pmatrix} \\
 & = 5+ 3(1+y+y^2)(1+z+z^2) \text{ mod } \langle y^3-1,z^3-1 \rangle= 5,
\end{align*} 
and
\begin{align*}  \alpha(y,z) & =\det \begin{pmatrix}1&1&1&1&1&1&1&1\\0 & 1 &  0 & z & 1 & 1 & 0 &0\\
0 & 0 & 1 & 0 & yz^2 & 1 &1 &0 \\
0 & 0 & 0 & 1 & 0 & y^2z^2 &1 & 1 \\
1 & 0 & 0 & 0 & 1 &0 & y^2 & 1 \\ 
1 & 1 & 0 &0 &0 & 1 & 0 & z^2\\
y^2z & 1 & 1 & 0 & 0 & 0 & 1 &0 \\ 
0 & yz & 1 & 1 & 0 & 0 & 0 & 1 \end{pmatrix} \\
 & = 3-y^2+y^2z^2-y^2z+yz^2+2yz+3z^2 \text{ mod } \langle y^3-1,z^3-1 \rangle.
\end{align*}
Using the 3 and $2yz$ (two terms with coefficients  differing by 1) we have
\begin{align*} \beta(y,z) & =(1-y^2z^2)\alpha(y,z)\\ & = 1+y-y^2-z+3z^2+yz-4y^2z+2yz^2-2y^2z^2  \;\; \bmod \langle y^3-1,z^3-1 \rangle,
\end{align*}
and  
$$\sum_{j=1}^8 \alpha(\vec{y}_j)=\sum_{i=0}^2\sum_{j=0}^2 \alpha(\om^i,\om^j)-\alpha(1,1)= 9\cdot 3 - 8=19,\quad \om:=e^{2\pi i/3}, $$
and $\sum_{j=1}^8 \beta(\vec{y}_j)=9\cdot 1 -0=9,$
giving $ B= 5+19c +9b.  $
\end{proof}

\section{Trying to generate  all $B\equiv A\pmod q$ for other $A$}

We have shown that for any $\mathbb Z_{q-1}$ determinant $A$ with $\gcd(A,q)=1$ or $q^2\mid A$ we can achieve a
$GA(1,q)$ determinant  $AB^{q-1}$ with any $B\equiv A\pmod q$. 
Observe also that for $q$ even $F\mapsto -F$ and for $q$ odd $F\mapsto xF$ sends $A,B\mapsto -A,-B$.

We will use the following lemma  to try to construct all $B\equiv A\pmod q$ for other 
values of $A$.

\begin{lemma}\label{Shifty}
Suppose that $F(x;y_1,\ldots ,y_n)$
has $A$ and $B$ values $A_0=A_1F(1;1,\ldots,1)$ and $B_0,$ and $g(x)$ has $\mathbb Z_{q-1}$ determinant $C_0=C_1g(1)\neq 0$. Then
$$ G(x;y_1,\ldots,y_k)=F(x;y_1,\ldots ,y_k)g(x)+ h(y_1,\ldots ,y_k) (1+x+\cdots + x^{q-2}) \mod (x^{q-1}-1),$$
has
$$ A=A_1C_1\left(F(1;1,\dots ,1)g(1)+(q-1)h(1,\ldots ,1)\right), \quad B=B_0C_0+C_1\left(q\lambda_{0\cdots 0} -A_1h(1,\ldots ,1)\right), $$
where, writing $\alpha (y_1,\ldots ,y_k)$ for the determinant of the matrix for $B_0$ but with the first row all changed to 1's, $\lambda_{0...0}$ is the constant term in 
\begin{align*} \lambda(y_1,\ldots ,y_k) &: =\alpha (y_1,\ldots ,y_k)h(y_1,\ldots ,y_k) \mod \langle y_1^p-1,\ldots ,y_k^p-1 \rangle \\
 & = \sum_{0\leq i_1,\ldots ,i_k <p} \lambda_{i_1...i_k}y_1^{i_1}\cdots y_k^{i_k}.
\end{align*}
\end{lemma}
Here $A_1$ means the product of $F(x;1,\ldots ,1)$  over the $(q-1)$st roots of unity excluding 1 (allowing  $F(1;1\ldots ,1)=0$).
\begin{proof}[Proof of Lemma \ref{Shifty}]  Observe that  for $g(x)=\sum_{i=0}^{q-2}c_ix^i$, the coefficient of $x^j$ in $F(x;y_1,\ldots ,y_k)g(x)$ will be the 
$\sum_{i=0}^{q-2} f_i(y_1,\ldots,y_k)c_{j-i \pmod {q-1}}.$
 In particular the matrix, $B_1$ say, for the corresponding $B$ will be the matrix for $F(x;y_1,\ldots ,y_k)$ for $B_0$,  multiplied on the right by the circulant determinant, $K$ say,  for $g(x)$. In particular $B_1$ will have determinant $B_0C_0$. 

Plainly, if we then add  $(1+x+\cdots +x^{q-2})h(y_1,\ldots ,y_k)$, the matrix for  $B$ will be the matrix for $B_1$ with $h(y_1,\ldots ,y_k)$ added to each entry  to the first row of  $B_1$ (likewise for the subsequent rows with
the variables $y_1,\ldots ,y_k$ running through the length $(q-1)$ cycle of monomials). 
Subtracting the last column from the first columns leaves $h(y_1,\ldots ,y_k)$ just in the last column and expanding down this last column will give the determinant of
the form  $\det(B_1)$
plus $h(y_1,\ldots ,y_k)\alpha'(y_1,\ldots ,y_k)$ for some $\alpha'(y_1,\ldots ,y_k)$, plus similar terms where the $y_1,\ldots ,y_k$ run through the cycle of monomials. It is not hard to see (do the same elimination) that $\alpha'(y_1,\ldots ,y_k)$ will be the determinant of the matrix obtained by replacing the first row of $B_1$ by all 1's.  
Notice, multiplying the matrix for $B_0$ with first row replaced by 1's (whose determinant is defined to be $\alpha(y_1,\ldots ,y_k)$) on the right by $K$ gives the matrix for $B_1$ with first row replaced by $g(1)$'s. Hence it's determinant, $g(1)\alpha'(y_1,\ldots ,y_k)$ is $\alpha(y_1,\ldots ,y_k)\det(K)$, 
and $\alpha'(y_1,\ldots ,y_k)=C_1\alpha(y_1,\ldots ,y_k)$ for any $s$.
Finally we need to sum the
$$\alpha(y_1,\ldots ,y_k)h(y_1,\ldots ,y_k)=\sum_{0\leq i_1,\ldots ,i_k<p} \lambda_{i_1...i_k}y_1^{i_1}\cdots y_k^{i_k}  \mod \langle y_1^p-1,\ldots ,y_k^p-1\rangle, $$ 
with the $y_1,\ldots ,y_k$ running through the $q-1$ cycle of monomials. If we sum a
$y_1^{i_1}\cdots y_k^{i_k}$, $0\leq i_1,\ldots ,i_k<p,$  over the sequence we get 
$(q-1)$ if $i_1=\cdots =i_k=0$ and $-1$ otherwise. Hence
$B=B_0C_0+ q\lambda_{0...0}- \lambda(1,\ldots ,1)$, where plainly $\lambda(1,\ldots ,1)=\alpha(1,\ldots ,1)h(1,\ldots,1)$. Notice that putting all the $y_i=1$ is the same 
as shifting $F(x;1,\ldots ,1)$ by  $m(1+x+\cdots +x^{q-1})$  say, reducing to the familiar case of shifting $\mathbb Z_{q-1}$ determinants. That is, $B_0=A_0$ and $B_0\mapsto B_0+m(q-1)\alpha(1,\ldots ,1)$ is the same as $A_0\mapsto A_1(F(1;1,\ldots ,1)+m(q-1))$, giving $\alpha(1,\ldots ,1)$.
\end{proof}
We will use this lemma as follows.

\noindent
\subsection{The Procedure}\label{procedure}
Suppose that $s$ in $\mathbb N$ has $\gcd(s,q-1)=1$ and that $A_0$ is a $\mathbb Z_{q-1}$ determinant. We try to construct all $B\equiv A\pmod q$ for all $A\equiv sA_0\mod {A_1(q-1)}$. Theorem \ref{coprime} and Theorem \ref{square} are special  cases of this process.

\vskip1ex
 Take a stem polynomial $f(x)$ in $\mathbb Z[x]$ achieving the $\mathbb Z_{q-1}$ determinant $A_0=A_1f(1)$.
We form a $F(x;y_1,\ldots, y_k)$ with $F(x;1,\ldots ,1)=f(x)$ and calculate 
$$\alpha(y_1,\ldots , y_k)=\sum_{0\leq l_1,\ldots ,l_n<p} \alpha_{l_1..l_k}y_1^{l_1}\cdots y_{k}^{l_k} \mod \langle y_1^p-1,\ldots y_k^p-1\rangle. $$
Our goal is to find a $t(y_1,\ldots ,y_k)$ such that $t(1,\ldots ,1)=0$ and
$$t(y_1,\ldots ,y_k)\alpha(y_1,\ldots ,y_k) \mod \langle y_1^p-1,\ldots y_k^p-1\rangle $$
has constant term 1. Then taking
$$G(x;y_1,\ldots ,y_k)=F(x;y_1,\ldots ,y_k)(1+\cdots +x^{s-1}) + \left( n+m t(y_1,\ldots ,y_k)\right) (1+\cdots +x^{q-2}),$$
reduced mod ${(x^{q-1}-1})$ if necessary, will give
$$ A=A_1\left( sf(1)+n(q-1)\right), \quad B=sB_0+n\left(q\alpha_{0...0}-A_1\right) + mq, $$
and for any $A$ of this form we can obtain any $B\equiv A \pmod q$ with a suitable 
choice of $m$.

If, by chance, $\alpha(y_1,\ldots ,y_k)$ has two coefficients differing by 1, 
$\alpha_{i_1...i_k}-\alpha_{i_1'...i_k'}=1$ say, then we can just take
$$ t(y_1,\ldots ,y_k)=y_1^{p-i_1}\cdots y_k^{p-i_k}-y_1^{p-i_1'}\cdots y_k^{p-i_k'}. $$
Otherwise we consider
$$ \beta (y_1,\ldots ,y_k)=(y_1-1)\alpha(y_1,\ldots ,y_k)=\sum_{0\leq i_1,\ldots ,i_k<p}\beta_{i_1...i_k}y_1^{i_1}\cdots y_k^{i_k} \mod \langle y_1^p-1,\ldots ,y_k^p-1\rangle. $$
Instead of $y_1-1$ we could take here any polynomial equaling 0 at $y_1=\cdots=y_k=1$.
If the gcd of the $\beta_{i_1...i_k}$ is 1 then we can find  a sum of integer multiples $\lambda_{i_1...i_k} \beta_{i_1...i_k}$ giving 1 and take 
$$ t(y_1,\ldots ,y_k)= (y_1-1)\sum_{0\leq i_1,\ldots ,i_k<p}\lambda_{i_1...i_k}y_1^{p-i_1}\cdots y_k^{p-i_k} \mod \langle y_1^p-1,\ldots ,y_k^p-1\rangle. $$
If the gcd is greater than 1 then we try a new $F(x;y_1,\ldots ,y_k)$.
Notice, if  the gcd of the $\beta_{i_1...i_k}$ is $\beta>1$ then this construction still gives 
\be \label{beta} A=A_1\left( sf(1)+n(q-1)\right), \quad B=sB_0+n\left(q\alpha_{0...0}-A_1\right) + m\beta q. \ee

Plainly, if we can find a suitable $t(y_1,\ldots ,y_k)$ for $A_0$ (and so obtain all
$A=A_0$ and $B\equiv A_0 \mod q$) we can do the same for $A=sA_0$, for any  $\gcd(s,q-1)=1$. In fact, if $p_1^{\alpha}\parallel q-1$ for some prime $p_1$  and $p_1^{\alpha}\mid f(1)$
then we obtain the result for $A=A_0t$ with any $\gcd(t,(q-1)/p_1^{\alpha})=1. $
To see this take $n=n_1 f(1)/p_1^{\alpha}$ and solve $s+n_1 (q-1)/p_1^{\alpha}=t$ with 
$\gcd(s,q-1)=1$  and $s>0$ (for example, by solving  $s\equiv t \pmod{(q-1)/p_1^{\alpha}}$ and $s\equiv 1\pmod {p_1}$ with $s>0$).

\noindent
\subsection{Special case:  $q-1=2\ell$ for some prime $\ell$. } 

From Theorem 5.1 and Theorem 5.2 we obtain 
all $B\equiv A\pmod q$ for any $A$ with $\gcd(A,2\ell)=1$ or $4\ell^2\mid A$. This just leaves
 the $\mathbb Z_{2\ell}$ determinants $A=\ell^2n,$ $n$ odd and the $A=4n$, $\ell\nmid n$. From the above comments 
we just need to show the existence of a suitable $t(y_1,\ldots ,y_k)$ for the two values
$A_0=\ell^2$ achieved with $\ell\mid f(1)$, for example (see \cite[Lemma 5]{Laquer}) take $f(x)= (x^{\ell+1}-1)/(x-1)-x $, and $A_0=4$ achieved with $2\mid f(1)$, for example take $f(x)=1+x^2$. In the next section we illustrate this for $GA(1,27)$.

\section{The group $GA(1,27)$}\label{q=27}
In this case we know that the integer group determinants take the form
\be \label{27} D=AB^{26}, \quad B\equiv A \bmod 27,\;\; \ee
where $A$ is a $\mathbb Z_{26}$ integer determinant; that is,  by \cite{Laquer},
\be \label{Z26}  2\mid A\Rightarrow 2^2\mid A,\quad 13\mid A\Rightarrow 13^2\mid A. \ee
We prove that these conditions are sufficient as well as necessary.

\begin{theorem}
The integer group determinants for $GA(1,27)$ are exactly the integers of the form \eqref{27}
with $A$ satisfying \eqref{Z26}.

\end{theorem}

\begin{proof} Theorems \ref{coprime} and \ref{square} prove this for the $A=m$ with $\gcd(m,26)=1$ and  the $A=26^2m$ with any $m\in \mathbb Z$. This leaves the $A=2^2m$ with $\gcd(m,13)=1$, and the $A=13^2m$ with $m$ odd.

\subsection{Achieving the remaining values} We use $y,z,w$ for $y_1,y_2,y_3.$ We  take $r$ to be a root of $f(x)=x^3-x+1$ over $\mathbb F_3$,
and calculate the exponents on $y,z,w$ corresponding to $y$ in the  $j$th row, by 
$$ M^{j-1}\begin{pmatrix} 1 \\ 0 \\ 0 \end{pmatrix} \bmod 3, \;\; j=1,\ldots,26, \quad M:=\left[ \begin{matrix} 0 & 0 & -1 \\ 1 & 0 & 1 \\ 0 & 1 & 0\end{matrix}\right]. $$
Hence as we go down the rows the variable $y$ runs through the sequence
\begin{align*} y,z,w, y^2z,z^2w, y^2zw^2,yzw,y^2z^2w,y^2w^2,yz,zw,y^2zw,y^2w,&\\
y^2,z^2,w^2,yz^2,zw^2,yz^2w,y^2z^2w^2,yzw^2,yw,y^2z^2,z^2w^2,yz^2w^2,yw^2, &
\end{align*}
and $z$ and $w$ the circulant shifts of this sequence by one or two places respectively.

As observed in the previous section we just have to deal with $A_0=2^2$, $2\mid f(1)$ and $A=13^2$, $13\mid f(1)$ and show the existence of a suitable $t(y,z,w)$.

\vskip1ex
\noindent
{\bf Case (i): $A_0=2^2$}. We took  $f(x)=1+x^2$ and $F(x;y,z,w)=1+(1-y)x+x^2$.
The $\alpha(y,z,w)$ did not have two coefficients differing by 1,
 so we calculated  $\beta(y,z,w)$
\begin{align*}&=  594 - 290 y - 304 y^2 - 163 z - 176 y z + 339 y^2 z - 192 z^2 + 
 100 y z^2 + 92 y^2 z^2 \\ & + 
 w (-483 + 126 y + 357 y^2 + 317 z + 97 y z - 414 y^2 z - 240 z^2 + 
    513 y z^2 - 273 y^2 z^2)\\ & + 
 w^2 (88 - 118 y + 30 y^2 - 265 z + 124 y z + 141 y^2 z + 32 z^2 + 
    3 y z^2 - 35 y^2 z^2)
\end{align*}
Using the $3yz^2w^2$  and $32z^2w^2$ we take $t(y,z,w)=(y-1)(11y^2zw-zw)$ and for $\gcd(s,26)=1$ get
$$ A=4(s+13n),\quad B=-1670s-13691n+27m. $$
Taking $s=1,15,3,17,5,19$ and $F\mapsto Fx$ give all the $A=4t$, $13\nmid t$.

\vskip1ex
\noindent
{\bf Case (ii): $A_0=13^2$}.
For $A_0=13^2$ we use 
$F(x;y,z,w)=1+(y-z)x+x^2+\cdots +x^{13}$. This produces $\beta(y,z,w)=$
\begin{align*}
 & -7947 + 2923 y + 5024 y^2 + 6394 z + 5132 y z - 11526 y^2 z + 
 3644 z^2 - 5820 y z^2 + 2176 y^2 z^2 \\ & + 
 w (2956 + 761 y - 3717 y^2 - 2457 z + 2379 y z + 78 y^2 z - 
    9319 z^2 + 4659 y z^2 + 4660 y^2 z^2) \\ & +
w^2 (-7343 - 1885 y + 9228 y^2 + 1098 z + 6361 y z - 7459 y^2 z + 
    426 z^2 - 1355 y z^2 + 929 y^2 z^2),
\end{align*}
and using the $78y^2zw$ and $761yw$ we take $t(y,z,w)=(y-1)yw^2(361z^2-37y)$,
and $G(x;y,z,w)=F(x;y,z,w)+(n + m t(y,z,w))(1+x+\cdots +x^{25})$, gives $A=13^2(1+2n)$ and $B=3436+ 29525n+27m$. 
\end{proof}

\section{GA(1,16)}\label{q=16}

Recall that the $GA(1,16)$ determinants take the form
$$D= AB^{15},\quad B\equiv A \bmod {16}, $$
where $A$ is a $\mathbb Z_{15}$ determinant. Unlike the previous groups we shall show that there are $A$ for which we can not achieve all $B\equiv A\pmod q$.
A complete description of the integer group determinants for $GA(1,16)$ becomes quite complicated, so we have not tried to obtain it here.

\begin{theorem} Suppose that $A=3^2\cdot 2^5$ or $5^2\cdot 2^5,$ then any $B=16k$
must have $k$ even. In particular we cannot have $B=16$.
\end{theorem}
Note, $3^2\cdot 2^5$ and $5^2\cdot 2^5$  are $\mathbb Z_{15}$ determinants,
for example take $f(x)=(1+x)(1+x+x^4)$ and $(1+x)(1+x+x^2-x^4+x^7+x^8+x^9)$ respectively

\begin{proof} We use $X;w,x,y,z$ for the variables $x;y_1,y_2,y_3,y_4$. We  take $r$ to be a root of $x^4+x+1 \pmod 2$  giving the cycle of monomials
\be \label{variableorder} w,x,y,z,wx,xy,yz,wxz,wy,xz,wxy,xyz,wxyz,wyz,wz, \ee
Suppose that
$$ F(X;x,y,z,w)=\sum_{n=0}^{14} f_n(w,x,y,z)X^j $$
where, ordering the monomials according to \eqref{variableorder},
\begin{align*}  f_n(w,x,y,z)=  c_n + c_{n0}(1-w)+c_{n1}(1-x)+c_{n2}(1-y)+\cdots 
+ c_{n14}(1-wz),   \quad c_{nj}\in\mathbb Z
\end{align*}
and 
$$ F(X;1,1,1,1)=\sum_{n=0}^{14} c_nX^n =:f(X). $$
As in Section \ref{alternative}, writing $S$ for the diagonal matrix with
\begin{align*}  (1-w),(1-x),(1-y),(1-z),(1-wx),(1-xy),(1-yz),(1-wxz), & \\(1-wy),(1-xz), (1-wxy),(1-xyz),(1-wxyz),(1-wyz),(1-wz) &\end{align*}
down the diagonal, and $P$ for the $15\times 15$ circulant shift matrix \eqref{defP},
the matrix for $B$
takes the form
$$ \sum_{n=0}^{14} c_nP^n + \sum_{n=0}^{14}\sum_{j=0}^{14}  c_{nj} P^j S P^{n-j}. $$
Putting $w=x=y=z=1$ we get the circulant determinant $\sum_{n=0}^{14} c_nP^n$ for $f(X)$ whose determinant is $A$, any other choice of $\pm 1$ gives determinant $B$. In particular, taking $w=x=y=z=-1$ we get
\be \label{Bform}  B= \det{\mathcal B}, \quad\quad {\mathcal B}= \sum_{n=0}^{14} c_nP^n + 2\sum_{n=0}^{14}\sum_{j=0}^{14}  c_{nj} P^j S_0 P^{n-j}, \ee
where $S_0$ is the diagonal matrix with diagonal entries
$$ d_0,d_1,\ldots ,d_{14}:= 1,1,1,1,0,0,0,1,0,0,1,1,0,1,0.$$

Suppose that $A=3^2\cdot 2^5$ or $3^2\cdot 2^5$ is the $\mathbb Z_{15}$ determinant of $f(X)$. For a primitive 15th root of unity $\om$ write
$A=f(1)N_3N_3N_{15}$ with integers
$$
N_3=f(\om^5)f(\om^{10}),\quad N_5=f(\om^3)f(\om^6)f(\om^9)f(\om^{12}),$$
and
$$ N_{15}=f(\om)f(\om^2)f(\om^{4})f(\om^8)f(\om^{-1})f(\om^{-2})f(\om^{-4})f(\om^{-8}). $$
Recall that 2 stays prime in the 3rd and 5th roots of unity but factors into two primes of norm $2^4$ in the 15th roots of unity. From \cite[Theorem 4]{Z15} we see that 
$(f(1),N_3,N_5,N_{15})=(2\cdot 3, 3, 1, 2^4)$ or $(2\cdot 5,1,5,2^4).$
 There are two extensions of the 2-adic absolute 
value to the 15th roots of unity (writing $\Phi_{15}(x)=(x^4+x+1)(x^4+x^3+1)-2x(x^3+1)^2$ we see that 
$ |\om^4+\om+1|_2,|\om^{-4}+\om^{-2}+1|_2=2^{-1},1$ or $1,2^{-1}$). 
We use $|\;\;|_2$ to denote the valuation with $|f(\om)|_2<1$. Since $g(\om^{2^j})\equiv g(\om)^{2^j}\pmod 2$ we see that the conjugates of $\om$ over the $p$-adics $\mathbb Z_2$ will be $\om^2$,$\om^4$ and $\om^8$ and
$$ |f(\om)|_2=|f(\om^2)|_2= |f(\om^4)|_2=|f(\om^8)|_2=2^{-1}.$$

Take $V$ to be the matrix whose $(j+1)$st column has entries $1,\om^j,\om^{2j},\ldots ,\om^{14j}$ for $j=0,\ldots 14,$ and $U$ the matrix whose $(i+1)$st row has entries
$1,\om^{-i},\om^{-2i},\ldots ,\om^{-14i}$ for $i=0,\ldots 14$. Plainly $V^{-1}=\frac{1}{15}U$. Observing that $PV$ multiplies the $(j+1)$st column of $V$ by $\om^j$ we get $V^{-1}PV$ is the diagonal matrix with diagonal entries 
$1,\om,\om^2,\ldots ,\om^{14}$ and $V^{-1}\left( \sum_{n=0}^{14} c_nP^n\right) V$ is the diagonal matrix with $f(1),f(\om),\ldots ,f(\om^{14})$ down the diagonal.
Similarly the entry in the $(r+1)$st row and $(\ell +1)$st column of $V^{-1}SV$ will be
$$ E_{r\ell} = \frac{1}{15} \sum_{i=0}^{14} d_i \om^{(\ell -r)i}, \quad 0\leq r,\ell \leq 14, $$
and the $(r+1)$st row and $(\ell+1)$st column entry of ${\mathcal B}'=V^{-1}{\mathcal B}V$
will be
$$ 2\sum_{n=0}^{14} \sum_{j=0}^{14} c_{nj} \om^{rj}E_{r\ell } \om^{(n-j)\ell},$$
plus $f(\om^r)$ on the diagonal  $r=\ell$.
Notice that we can factor out a 2 from the rows corresponding to $f(\om^j)$, $j=0,1,2,4$ and 8 and still be left with all entries $|b_{ij}|_2\leq 1$. So 
$|B|_2=|\det{\mathcal B}' |_2 \leq 2^{-5},$ and the integer $B$ must be  a multiple of 32.
\end{proof}

\end{document}